\documentclass{amsart}
\usepackage{graphicx}
\usepackage{amscd}
\usepackage{amsmath}
\usepackage{amsfonts}
\usepackage{amssymb}
\newtheorem{theorem}{Theorem}
\theoremstyle{plain}

\newtheorem{example}{Example}[section]

\newtheorem{lemma}{Lemma}[section]

\newtheorem{remark}{Remark}[section]

\numberwithin{equation}{section}

\begin{document}
\title[Semilinear Equations at a Critical Exponent]{The Structure of the Solutions to Semilinear Equations at a Critical Exponent}
\author{Allan L. Edelson}
\address{Department of Mathematics\\
University of California\\
Davis, CA 95616}
\email{aedelson@math.ucdavis.edu}
\dedicatory{Dedicated to Vic Shapiro on the occasion of his retirement }

\begin{abstract}
This paper is concerned with the structure of the solutions to subcritical
elliptic equations related to the Matukuma equation. In certain cases the
complete structure of the solution set is known, and is comparable to that of
the original Matukuma equation. Here we derive sufficient conditions for a
more complicated solution set consisting of; (i) crossing solutions for small
initial conditions and large initial conditions; (ii) at least one open
interval of slowly decaying solutions; and (iii) at least two rapidly decaying
solutions. As a consequence we obtain multiplicity results for rapidly
decaying, or minimal solutions.
\end{abstract}\maketitle

\section{Introduction}

This article is concerned with the structure of the solutions to the
semilinear elliptic equation%

\begin{equation}
-\Delta u=f(|x|)u^{p},x\in{\mathbf{R}}^{n},\quad n\geq3,\quad\frac{n+2}%
{n-2}>p>1. \tag{1.1}%
\end{equation}
It is assumed that $f$ is positive and has asymptotic behavior $f(|x|)\sim
|x|^{l},-2<l<0$ as $x\rightarrow\infty$. Among the important examples of
related equations are the Matukuma equation%

\[
-\Delta u = {\frac{1}{1+|x|^{2}}} u^{p}, \quad x\in{\mathbf{R}} ^{3},
\]
and the scalar curvature equation%

\[
-\Delta u = f(x)u^{\frac{n+2}{n-2} }.
\]

\noindent There has been a substantial body of research into the question of
the existence or non existence of positive solutions of (1.1) decaying to zero
as $|x|\rightarrow\infty$. Much of this research has focussed on the range of
$p$ values between the Sobolev critical exponent $p^{\ast}=\frac{n+2}{n-2}$,
and the critical exponent%

\begin{equation}
p_{\ast}=\frac{n+2+2l}{n-2}. \tag{1.2}%
\end{equation}
The exponent $p_{\ast}$ is critical with respect to the compact imbedding of
weighted Sobolev spaces%

\[
H^{1}_{r}({\mathbf{R}}^{n})\rightarrow L^{q}({\mathbf{R}}^{n},f(|x|)dx),
0<q<p_{*},
\]
(see [Eg]). For $p>p_{*}$ variational methods are sufficient to prove the
existence of positive, decaying solutions (see [Eg], [N-S]). For $p=p_{*}$
those methods fail due to non compactness, and other techniques are required.
This article will be primarily concerned with the structure of the set of
solutions when $p=p_{*}$.

\bigskip

\noindent When $f$ decays rapidly to zero at infinity, e.g. $l<-2$, existence
and nonexistence results for (1.1) were obtained in [N]. References for the
case $-2<l<0$, include [K-N-Y], [N-Y], [N-S], [Y-Y1] and [Y-Y2]. For
$0<f(0)<\infty$ and $p\neq p_{\ast}$, the existence or nonexistence of ground
states is determined by the decay rate $l$. At the critical exponent,
$p=p_{\ast}$, the existence of positive, decaying solutions depends on higher
order terms in the asymptotic expansion of $f(r)$. In particular, the function%

\begin{equation}
h(r)=r(r^{-l}f(r))_{r}, \tag{1.3}%
\end{equation}
has been useful in obtaining more precise information about the structure of
the set of solutions when $0<f(0)\leq\infty$. The following examples
illustrate the role played by $p$ and $h(r)$ in determining the existence or
nonexistence of positive, entire solutions. We assume that $f$ is radially
symmetric and so the problem reduces to an ordinary differential equation.
Since we are interested in nonnegative solutions we define $u^{+}=max\{u,0\}$
and let $u(r;\alpha)$ be the (unique) solution of the initial value problem%

\begin{equation}
u^{\prime\prime}+\frac{n-1}{r}u^{\prime}+f(r)(u^{+})^{p}=0,\quad
u(0;\alpha)=\alpha,\quad u^{\prime}(0;\alpha)=0. \tag{1.4}%
\end{equation}

\begin{example}
$f(r)=r^{l},-2<l<0$.

\begin{itemize}
\item  If $1<p<p_{\ast}$, then $u(r;\alpha)$ has a finite zero for every
$\alpha>0$.

\item  If $p_{\ast}\leq p$ then $u(r;\alpha)$ is a positive, decaying solution
for every $\alpha>0$.
\end{itemize}

In this case $h(r)=0$, for $r\geq0$.
\end{example}

\begin{example}
$f(r)=O(r^{\sigma})$ as $r\rightarrow0$, $\sigma>-2$.

\begin{itemize}
\item  If $p=p_{\ast}$, $f(r)\geq cr^{l}$ for large $r$, and $\int_{0}%
^{R}h(r)r^{n+l-1}dr\leq0$, for all $R>0$, then for every $\alpha>0$,
$u(r;\alpha)$ is a positive, entire solution satisfying $lim_{r\rightarrow
\infty}u(r;\alpha)=0$.

\item  If $p=p_{\ast}$, $h(r)\not \equiv0$, and $\int_{0}^{R}h(r)r^{n+l-1}%
dr\geq0$ for all $R>0$, then $u(r;\alpha)$ has a finite zero for all
$\alpha>0$ .
\end{itemize}
\end{example}

\begin{example}
Let $\frac{n+l}{n-2}<p<p_{\ast}$, and
\[
f(r)=\frac{(l+2)(n-2)}{(p-1)^{2}}\left\{  \left(  p-\frac{n+l}{n-2}\right)
+\left(  \frac{l+2p}{n-2}\right)  \left(  \frac{1}{1+r^{2}}\right)  \right\}
(1+r^{2})^{l/2}.
\]

\begin{itemize}
\item $u(r;\alpha)$ has a finite zero for every sufficiently large $\alpha>0$,

\item $u(r;1)=(1+r^{2})^{-\frac{l+2}{2(p-1)}}$ is a positive, decaying solution,

\item $u(r;\alpha)$ has a finite zero for every sufficiently small $\alpha>0$.
\end{itemize}
\end{example}

Example I is contained in Proposition 4.5 of [N-Y]. The results of Example II
can be found in Theorems 9.2 and 9.3 of [K-N-Y], where it is assumed that
$\int_{0}^{R}h(r)r^{n+l-1}dr$ does not change sign. Example III shows the
complexity of the solutions when $h(r)$ and $\int_{0}^{R}h(r)r^{n+l-1}dr$ are
allowed to change sign (see also [N: Example 5.2], [N-S]). An important
structure theorem for equation (1.1) is given in [K-Y-Y], where it is assumed
that this integral changes sign exactly once. \smallskip

\noindent It is known ([Y-Y1]) that under our hypotheses all solutions of
(1.4) with $\alpha>0$, are of one of the following types: \smallskip

\begin{enumerate}
\item [(1)]$u(r;\alpha)$ is a crossing solution, i.e. has a positive zero,

\item[(2)] $u(r;\alpha)$ is a slowly decaying solution, i.e.%
\[
lim_{r\rightarrow\infty}r^{n-2}u(r;\alpha)=\infty,
\]

\item[(3)] $u(r;\alpha)$ is a rapidly decaying solution, i.e.%
\[
\overline{lim}_{r\rightarrow\infty}r^{n-2}u(r;\alpha)<\infty.
\]
\end{enumerate}

In sections 3 and 4 we will give criteria for existence and non existence of
positive, decaying solutions when both $h(r)$ and $\int_{0}^{R}h(r)r^{n+l-1}%
dr$ are allowed to change sign. Then in section 5 we prove that for
$p=p_{\ast}$, and under suitable hypotheses on $h(r)$, the set of solutions of
(1.4) consists of:

\begin{enumerate}
\item [(1)]crossing solutions for sufficiently small or sufficiently large
$\alpha$;

\item[(2)] at least one open interval of $\alpha$ values for which
$u(r;\alpha)$ is slowly decaying;

\item[(3)] at least two rapidly decaying solutions.
\end{enumerate}

\noindent Thus under these hypotheses a structure comparable to that of
Example III can be achieved for $p=p_{*}$. A consequence of this result is
that we prove the existence of multiple minimal, i.e. rapidly decaying
solutions. \bigskip

\section{Main Results}

The following hypotheses will be assumed throughout this paper.

\noindent$(f_{1})\quad f\in C(0,\infty),\quad f(r)>0\quad$ for $r>0,$

\noindent$(f_{2})\quad f(r)=O(r^{l})\quad$as $r\rightarrow\infty,-2<l<0,$

\noindent$(f_{2}^{\prime})\quad f(r)=O(r^{\sigma})\quad$as $r\rightarrow
0,-2<\sigma,$ \smallskip

\noindent We will assume additional hypotheses chosen from the following:

\noindent$(f_{3})\quad$ there exist positive numbers $\delta_{1},r_{2},\beta$
such that%

\[
\qquad h(r)<-\delta_{1}r^{-\beta},\quad\text{for }r>r_{2},
\]

\noindent$(f_{4})\quad h(r)=O(r^{\gamma}),\quad\gamma>0,\quad$as $r\rightarrow0.$

\noindent$(f_{5})\quad\int_{0}^{\infty}h(r)r^{n+l-1}dr<0$.

\noindent$(f_{6})\quad$ there exist positive numbers $\delta^{\prime}%
_{1},r^{\prime}_{2},\beta$ such that%

\[
\qquad|h(r)|<\delta_{1}^{\prime}r^{-\beta},\quad\text{for }r>r_{2}^{\prime},
\]

\noindent$(f_{7})\quad\int_{0}^{\infty}h(r)r^{n+l-1}dr>0,$

\noindent$(f_{8})\quad0<r_{3}=sup\{r>0:\int_{o}^{r} h(s)s^{n+l-1}ds \le0 \}.$

\noindent$(f_{9})\quad\gamma(2+l)>(\sigma-l)(n+l).$

\bigskip

\noindent For $p\ge p_{*}$ the following two theorems give sufficient
conditions for the existence of positive, entire solutions when $\alpha$ is
sufficiently small.

\begin{theorem}
Assume that ($f_{3}$) and ($f_{4}$) are satisfied, $p\geq p_{\ast}$ and that
$0<\beta<n+l$. Then there exists an $\alpha_{0}>0$ such that for $0<\alpha
\leq\alpha_{0}$, the solution $u(r;\alpha)$ is a positive, entire solution
which also satisfies $lim_{r\rightarrow\infty}u(r;\alpha)=0$.
\end{theorem}

\begin{example}
The function%
\begin{align*}
&  f(r)=(c_{1}+c_{2}r^{2})^{\tfrac{\gamma}{2}}\left(  c_{3}+c_{4}r^{2}\right)
^{\tfrac{\nu}{2}}\\
-2 &  <\gamma+\nu<0\text{, }c_{1}c_{4}\gamma+c_{2}c_{3}\nu>0\text{, }c_{i}>0
\end{align*}
satisfies $(f_{1})-(f_{4})$.
\end{example}

\begin{example}
The function%
\[
f(r)=(1+(1+r^{2})^{-\frac{1}{4}})(1+r^{2})^{-\frac{1}{4}}%
\]
satisfies $\left(  f_{1}\right)  -\left(  f_{4}\right)  $. For $n=3$ we have
$l=-\frac{1}{2},p_{\ast}=4,\beta=\frac{1}{2}$, and $0<\beta<n+l$.
\end{example}

\noindent This example is not covered by [K-N-Y], Theorem 9.2, since $\int
_{0}^{R}h(r)r^{n+l}dr>0$ for small $R$.

\begin{example}
The function%
\[
f(r)=\left(  \frac{9}{8}+(1+r^{2})^{-1}\right)  (1+r^{2})^{-\frac{3}{4}}%
\]
satisfies ($f_{1}$)--($f_{4}$). For $n=3$ we have $l=-\frac{3}{2},p_{\ast
}=2,\beta=2$, and $\beta>n+l$. Since $\int_{0}^{R}h(r)r^{n+l-1}dr>0$ for all
$R>0$, it follows from Theorem~9.3 of [K-N-Y] that $u(r;\alpha)$ has a
positive zero for all $\alpha>0$.
\end{example}

\noindent The next result gives a variant of Theorem 1 in which the condition
$\beta<n+l$ is replaced by an integral growth condition.

\begin{theorem}
Assume that ($f_{3}$), ($f_{4}$) and ($f_{5}$) are satisfied, and $p\geq
p_{\ast}$. Then there exists an $\alpha_{0}>0$ such that for $0<\alpha
\leq\alpha_{0}$, the solution $u(r;\alpha)$ is a positive, entire solution
which also satisfies $lim_{r\rightarrow\infty}u(r;\alpha)=0$.
\end{theorem}

\noindent The next two theorems give sufficient conditions for the non
existence of positive solutions when $\alpha$ is sufficiently large, or
sufficiently small. \bigskip

\begin{theorem}
Assume that $(f_{4})$, $(f_{6})$ and $(f_{9})$ are satisfied, and $p\geq
p_{\ast}$. Then there exist positive numbers $\gamma_{\ast},\alpha_{1}$, such
that for $\alpha>\alpha_{1}$ and $0<\gamma<\gamma_{\ast}$ the solution
$u(r;\alpha)$ has a positive zero.
\end{theorem}

\begin{remark}
For $p<p_{\ast}$ this result is proved in [N-Y, Theorem 2], without the
additional hypothesis $\gamma<\gamma_{\ast}$.
\end{remark}

\begin{theorem}
Assume that $(f_{4})$ and $(f_{7})$ are satisfied, $p\geq p_{\ast}$ and that
there exists an $r_{2}$ such that $h(r)\geq0,$ for $r_{2}<r<\infty$. Then
there exists an $\alpha_{0}>0$ such that for $\alpha\in(0,\alpha_{0})$ the
solution $u(r;\alpha)$ has a positive zero.
\end{theorem}

\begin{remark}
A stronger version of this result is found in [Y-Y2, Theorem 3].
\end{remark}

\begin{theorem}
Let $p=p_{\ast}$ and let $k(r)$ be a continuous function satisfying $(f_{7})$,
$(f_{9})$ and also the following:

\noindent there exist positive numbers $0<a<b<c$, $\alpha_{\ast}$, and
$r_{\ast}>b$ such that%
\begin{align}
(a)  &  \quad k(0)=k(a)=k(b)=k(c)=0,\quad k(r)=0,\text{ for }r>c,\tag{2.1}\\
(b)  &  \quad k(r)>0\text{ for }0<r<a,\text{ and }b<r<c,\quad k(r)<0\text{ for
}a<r<b,\nonumber\\
(c)  &  \quad k(r)=O(r^{\gamma}),\quad n+l-1<\gamma<\gamma_{\ast},\text{ as
}r\rightarrow0,\nonumber\\
(d)  &  \quad\phi(0;\alpha_{\ast})^{p_{\ast}+1}\int_{0}^{a}k(s)ds+\phi
(b;\alpha_{\ast})^{p_{\ast}+1}\int_{a}^{b}k(s)ds\nonumber\\
&  \quad+\phi(b;\alpha_{\ast})^{p_{\ast}+1}\int_{b}^{r_{\ast}}k(s)ds<-\delta
^{2},\nonumber\\
(e)  &  \quad\int_{0}^{c}k(r)dr>0,\nonumber
\end{align}
where $\gamma_{\ast}$ is given in theorem 3 and $\phi(r;\alpha)$ is the
solution of the equation%
\begin{equation}
-(r^{n-1}\phi_{r})_{r}=r^{n-1}r^{l}\phi^{p_{\ast}},\newline \phi
(0)=\alpha,\quad\phi_{r}(0)=0. \tag{2.2}%
\end{equation}
Then there exists a positive function $f\in C(0,\infty)$, with $f(r)=O(r^{l})$
as $r\rightarrow\infty,-2<l<0$, and positive numbers $\alpha_{0}<\alpha
_{1}<\alpha_{\ast}<\alpha_{2}<\alpha_{3}$ such that the solutions
$u(r;\alpha)$ of equation (1.4) are of the following type:%
\begin{align}
(a)  &  \quad0<\alpha<\alpha_{0}\Rightarrow u(r;\alpha)\text{ as a positive
zero,}\tag{2.3}\\
(b)  &  \quad u(r;\alpha_{\ast})\text{ as a positive, slowly decaying
solution,}\nonumber\\
(c)  &  \quad\alpha_{3}<\alpha<\infty\Rightarrow u(r;\alpha)\text{ as a
positive zero,}\nonumber\\
(d)  &  \quad u(r;\alpha_{1})\text{ and }u(r;\alpha_{2})\text{ are rapidly
decaying solutions}.\nonumber
\end{align}
\end{theorem}

\noindent It is not difficult to construct examples which satisfy the
hypotheses (2.1).

\section{Existence of positive, decaying solutions}

The basic existence theorem for solutions to the initial value problem can be
found in [N-Y]. In establishing Theorem~1 we will require several preliminary
lemmas. The first is contained in Proposition 4.1 of [N-Y] and is also found
in [N].

\begin{lemma}
There exists a unique solution $u(r;\alpha)\in C((0,\infty)\cap C^{2}%
(0,\infty))$ of (1.4). The solution satisfies the integral equation%
\begin{equation}
u(r)=\alpha-\frac{1}{n-2}\int_{0}^{r}\left\{  1-\left(  {\frac{s}{r}}\right)
^{n-2}\right\}  sf(s)(u^{+}(s))^{p}ds, \tag{3.1}%
\end{equation}
and also the following conditions;%
\begin{align}
(a)  &  \quad lim_{r\rightarrow\infty}r^{n-1}u_{r}(r)=0\newline (b)\quad
u_{r}(r)=-\int_{0}^{r}\left(  \frac{s}{r}\right)  ^{n-1}f(s)(u^{+}%
(s))^{p}ds\leq0,\quad\text{for }r>0,\tag{3.2}\\
(c)  &  \quad u\text{ is non-increasing on }[0,\infty).\nonumber
\end{align}
\end{lemma}

\begin{lemma}
([N-Y], Prop.4.3 and (6.8) ) For $R>0$ the solution $u(r;\alpha)$ satisfies
the Pohozaev type identity%
\begin{align}
&  \frac{n-2}{2}{R}^{n-1}u(R;\alpha)u^{\prime}(R;\alpha)+\frac{1}{2}%
R^{n}u^{\prime}(R;\alpha)^{2}+\frac{1}{p+1}R^{n}f(R)(u^{+}(R;\alpha
))^{p+1}\tag{3.3}\\
&  =\frac{1}{p+1}\int_{0}^{R}\left\{  -\frac{n-2}{2}(p-p_{\ast})r^{-l}%
f(r)+h(r)\right\}  r^{n+l-1}(u^{+}(r;\alpha))^{p+1}dr.\nonumber
\end{align}
\noindent Define $r_{\alpha}$ by
\end{lemma}%

\begin{equation}
r_{\alpha}=inf\left\{  r>0:u(r;\alpha)=\frac{\alpha}{2}\right\}  . \tag{3.4}%
\end{equation}

\begin{lemma}
$r_{\alpha}$ satisfies the following estimates:%
\begin{align}
r_{\alpha}  &  =O(\alpha^{\frac{1-p}{2+l}})\text{ as }\alpha\rightarrow
0,\tag{3.5a}\\
r_{\alpha}  &  =O(\alpha^{\frac{1-p}{2+\sigma}})\text{ as }\alpha
\rightarrow\infty. \tag{3.5b}%
\end{align}
\end{lemma}

\begin{proof}
By Proposition 4.1 of [N-Y] and ($f_{2}$) we have%
\begin{align*}
\frac{n-2}{2}\alpha &  =\int_{0}^{r_{\alpha}}\left\{  1-\left(  \frac
{s}{r_{\alpha}}\right)  ^{n-2}\right\}  sf(s)(u^{+}(s;\alpha))^{p}ds\\
&  \leq\alpha^{p}\int_{0}^{r_{\alpha}}sf(s)ds.\\
\frac{n-2}{2}\alpha^{1-p}  &  \leq\int_{0}^{r_{\alpha}}sf(s)ds.
\end{align*}
By Lemma 6.1 of [N-Y] $lim_{\alpha\rightarrow0}r_{\alpha}=\infty$, and
therefore%
\begin{align*}
\frac{n-2}{2}\alpha^{1-p}  &  \leq\int_{0}^{1}sf(s)ds+\int_{1}^{r_{\alpha}%
}sf(s)ds\\
&  \leq c+\int_{1}^{r_{\alpha}}s^{1+l}ds\leq cr_{\alpha}^{2+l},
\end{align*}
which implies%
\[
C\alpha^{\frac{1-p}{2+l}}\leq r_{\alpha}\text{ as }\alpha\rightarrow0.
\]
By a similar estimate we can show%
\[
C\alpha^{\frac{1-p}{2+l}}\geq r_{\alpha},\text{ as }\alpha\rightarrow0,
\]
and this proves 3.5a.

\noindent Consider the case $\alpha\rightarrow\infty$.%
\begin{align*}
\frac{n-2}{2}\alpha &  =\int_{0}^{r_{\alpha}}\left\{  1-\left(  \frac
{s}{r_{\alpha}}\right)  ^{n-2}\right\}  sf(s)(u(s;\alpha)^{+})^{p}ds\\
&  \geq\left(  \frac{\alpha}{2}\right)  ^{p}C\int_{0}^{r_{\alpha}}sf(s)ds,\\
\alpha^{1-p}  &  \geq\int_{0}^{r_{\alpha}}sf(s)ds,
\end{align*}
which implies $r_{\alpha}\rightarrow0$ as $\alpha\rightarrow\infty$. (3.5b)
can then be deduced.
\end{proof}

\begin{proof}
[Proof of Theorem 1]By $(f_{4})$ $h(r)>0$ on an interval $(0,\epsilon)$.
Define the quantities%
\begin{align}
r_{0}  &  =inf\{r>0:h(r)<0\},\tag{3.7}\\
r_{1}  &  =sup\{r>0:h(r)>0\},\nonumber\\
\delta_{2}  &  =\int_{0}^{r_{1}}|h(r)|r^{n+l-1}dr,\nonumber\\
k  &  =\left(  \frac{2^{-(p+1)}}{n+l-\beta}\right)  \left(  1-2^{-(n+l-\beta
)}\right)  >0.\nonumber
\end{align}
Then $r_{1}$ exists by $(f_{3})$, and by ($f_{4}$) $0<r_{0}\leq r_{1}<\infty$.
We apply the identity (3.3). By (3.5a) we may choose $\alpha_{0}$ so small
that%
\begin{equation}
r_{\alpha_{0}}^{n+l-\beta}\geq\frac{\delta_{2}}{k\delta_{1}}. \tag{3.8}%
\end{equation}
For some $\alpha\in(0,\alpha_{0}]$ assume that $u(r;\alpha)$ satisfies the
conditions%
\[
u(r;\alpha)>0,\quad0\leq r<R,\quad u(R;\alpha)=0.
\]
The left hand side of (3.3) is $>0$, because $u^{\prime}(R;\alpha)\neq0$.
Also, $R>r_{\alpha}$ by the definition of $r_{\alpha}$, and by Lemma (3.3) we
may assume $r_{1}<\frac{r_{\alpha}}{2}$. Then%
\begin{align}
\int_{0}^{R}h(r)r^{n+l-1}(u^{+}(r;\alpha))^{p+1}dr  &  =\tag{3.9}\\
\left\{  \int_{0}^{r_{1}}+\int_{r_{1}}^{\frac{r_{a}}{2}}+\int_{\frac{r_{a}}%
{2}}^{r_{\alpha}}+\int_{r_{\alpha}}^{R}\right\}  h(r)r^{n+l-1}(u^{+}%
(r;\alpha))^{p+1}dr  &  \leq\nonumber\\
\delta_{2}\alpha^{p+1}+\int_{\frac{r_{a}}{2}}^{r_{\alpha}}h(r)r^{n+l-1}%
(u^{+}(r;\alpha))^{p+1}dr  &  \leq\nonumber\\
\delta_{2}\alpha^{p+1}-\delta_{1}\int_{\frac{r_{a}}{2}}^{r_{\alpha}%
}r^{n+l-1-\beta}(u^{+}(r;\alpha))^{p+1}dr  &  \leq\nonumber\\
\delta_{2}\alpha^{p+1}-\delta_{1}\left(  \frac{\alpha}{2}\right)  ^{p+1}%
\int_{\frac{r_{a}}{2}}^{r_{\alpha}}r^{n+l-1-\beta}dr  &  \leq\nonumber\\
\alpha^{p+1}\left(  \delta_{2}-k\delta_{1}{r_{\alpha}}^{n+l-\beta}\right)   &
\leq0.\nonumber
\end{align}
This contradiction implies that $u(r;\alpha)$ is positive on $[0,\infty)$. The
condition%
\[
lim_{r\rightarrow\infty}u(r;\alpha)=0
\]
is a consequence of [N; Theorem 3.10].
\end{proof}

\begin{example}
If%
\[
f(r)=A_{0}r^{l}+o(r^{l}),\quad r\rightarrow\infty,
\]
then the hypotheses of Theorem 1 are satisfied.
\end{example}

\begin{proof}
It is easy to verify that in this case $\beta=1<n+l$.
\end{proof}

\noindent The next result, needed in the proof of Theorem 2, is a
generalization of Lemma 3.3.

\begin{lemma}
Define $r_{\alpha,k}$ by%
\[
r_{\alpha,k}=inf\left\{  r>0:u(r;\alpha)=\frac{\alpha}{k}\right\}  ,\quad
k>1.
\]
Then $r_{\alpha,k}$ has asymptotic behavior%
\[
lim_{\alpha\rightarrow0}r_{\alpha,k}=\infty.
\]
\end{lemma}

\begin{proof}
Following the proof of Lemma 3.3 we have%
\begin{align*}
u(r_{\alpha,k};\alpha)  &  =\frac{\alpha}{k}=\alpha-\frac{1}{n-2}\int
_{0}^{r_{\alpha,k}}\left\{  1-\left(  \frac{s}{r_{\alpha,k}}\right)
^{n-2}\right\}  sf(s)(u^{+}(s;\alpha))^{p}ds,\\
\alpha\left(  \frac{k-1}{k}\right)   &  =\frac{1}{n-2}\int_{0}^{r_{\alpha,k}%
}\left\{  1-\left(  \frac{s}{r_{\alpha,k}}\right)  \right\}  sf(s)(u^{+}%
(s,\alpha))^{p}ds\\
&  \leq\frac{1}{n-2}\int_{0}^{r_{\alpha,k}}sf(s)(u^{+}(s;\alpha))^{p}ds\\
&  \leq\left(  \frac{\alpha}{k}\right)  ^{p}\int_{0}^{r_{\alpha,k}}sf(s)ds,\\
\alpha^{1-p}  &  \leq C\int_{0}^{r_{\alpha,k}}sf(s)ds.
\end{align*}
Therefore%
\begin{align*}
\alpha^{1-p}  &  \leq C\int_{0}^{r_{\alpha,k}}sf(s)ds\\
&  =C\int_{0}^{r_{\alpha,k}}s^{l+1}ds\\
&  =Cr_{\alpha,k}^{l+2}.
\end{align*}
It follows that%
\[
\alpha^{\frac{1-p}{l+2}}\leq Cr_{\alpha,k}%
\]
and since $\frac{1-p}{l+2}<0$, we have the stated result.
\end{proof}

\noindent It is possible to obtain more precise estimates, analogous to
Lemma~3.3, however such estimates will not be needed.

\begin{proof}
[Proof of Theorem 2]Define%
\begin{align}
(a)  &  \quad r^{0}=sup\{r>0:\int_{0}^{r}h(s)s^{n+l-1}ds\geq0\},\tag{3.10}\\
(b)  &  \quad\Omega_{0}^{+}=\{r\in\lbrack0,r^{0}]:h(r)\geq0\},\nonumber\\
(c)  &  \quad\Omega_{0}^{-}=\{r\in\lbrack0,r^{0}]:h(r)\leq0\}.\nonumber
\end{align}
By ($f_{4}$) and ($f_{5}$) $r^{0}>0$. For $\epsilon\in(0,1)$ it is a
consequence of Lemma 3.5 that there exists an $\alpha_{0}(\epsilon)$ such that
for $\alpha\in(0,\alpha_{0}(\epsilon))$ the conditions%
\begin{align}
(i)  &  \quad(1-\epsilon)^{\frac{1}{p_{\ast}+1}}\alpha\leq u(r;\alpha
)\leq\alpha,\text{ for }0\leq r\leq r^{0},\tag{3.11}\\
(ii)  &  \quad r_{\alpha}>2r^{0},\nonumber\\
(iii)  &  \quad\int_{r^{0}}^{2r^{0}}h(r)r^{n+l-1}dr<2^{p_{\ast}+1}\epsilon
\int_{\Omega_{0}^{-}}h(r)r^{n+l-1}dr\nonumber
\end{align}
are satisfied. We will show that there exists an $\alpha_{0}$ such that
\[
\alpha\in(0,\alpha_{0})\Rightarrow\int_{0}^{R}h(r)r^{n+l-1}(u^{+}%
(r;\alpha))^{p_{\ast}+1}dr<0,\text{ for all }R>r^{0}.
\]
In fact, for $R>r^{0}$,
\begin{align*}
&  \int_{0}^{R}h(r)r^{n+l-1}(u^{+}(r;\alpha))^{p_{\ast}+1}dr\\
&  \leq\left(  \int_{\Omega_{0}^{-}}+\int_{\Omega_{0}^{+}}+\int_{r^{0}}%
^{R}\right)  h(r)r^{n+l-1}(u^{+}(r;\alpha))^{p_{\ast}+1}dr\\
&  \leq\alpha^{p_{\ast}+1}\int_{0}^{r^{0}}h(r)r^{n+l-1}dr-\epsilon
\alpha^{p^{\ast}+1}\int_{\Omega_{0}^{-}}h(r)r^{n+l-1}dr\\
&  +\int_{r^{0}}^{R}h(r)r^{n+l-1}dr+\int_{r^{0}}^{R}h(r)r^{n+l-1}%
(u^{+}(r;\alpha))^{p_{\ast}+1}dr\\
&  \leq-\epsilon\alpha^{p_{\ast}+1}\int_{\Omega_{0}^{-}}h(r)r^{n+l-1}%
dr+\int_{r^{0}}^{R}h(r)r^{n+l-1}(u^{+}(r;\alpha))^{p^{\ast}+1}dr.
\end{align*}
Then for $R>r^{0}$,%
\begin{align*}
&  -\epsilon\alpha^{p^{\ast}+1}\int_{\Omega_{0}^{-}}h(r)r^{n+l-1}%
dr+\int_{r^{0}}^{R}h(r)r^{n+l-1}(u^{+}(r;\alpha))^{p_{\ast}+1}dr\\
&  \leq-\epsilon\alpha^{p^{\ast}+1}\int_{\Omega_{0}^{-}}h(r)r^{n+l-1}%
dr+\int_{r^{0}}^{2r^{0}}h(r)r^{n+l-1}(u^{+}(r;\alpha))^{p_{\ast}+1}dr\\
&  \leq-\epsilon\alpha^{p_{\ast}+1}\int_{\Omega_{0}^{-}}h(r)r^{n+l-1}%
dr+\left(  \frac{\alpha}{2}\right)  ^{p_{\ast}+1}\int_{r^{0}}^{2r^{0}%
}h(r)r^{n+l-1}dr\\
&  \leq\alpha^{p_{\ast}+1}\left\{  -\epsilon\int_{\Omega_{0}^{-}}%
h(r)r^{n+l-1}dr+2^{-(p_{\ast}+1)}\int_{r^{0}}^{2r^{0}}h(r)r^{n+l-1}dr\right\}
<0,
\end{align*}
by (3.11iii).
\end{proof}

\noindent As in the proof of Theorem 1, an application of the Pohozaev
identity (3.3) leads to a contradiction, and once again the limit condition
follows from [N; Theorem 3.10].

\section{Nonexistence of positive, decaying solutions}

Our proofs of the nonexistence of positive, decaying solutions are based on
the Pohozaev identity. Define $w(r;\alpha)=r^{\frac{n-2}{2}}u(r;\alpha)$, and
let $u(r;\alpha)$ be a positive, entire solution which decays to zero as
$r\rightarrow\infty$. According to Theorem 2.2 of [K-N-Y] the following
variant of the Pohozaev identity applies.%

\begin{align}
&  -R^{2}w^{\prime\prime}(R;\alpha)w(R;\alpha)-Rw(R;\alpha)w^{\prime}%
(R;\alpha)+R^{2}w^{\prime}(R;\alpha)^{2}\tag{4.1}\\
&  -\left(  \frac{l+2}{n+l}\right)  R^{n}f(R)u(R;\alpha)^{p+1}\nonumber\\
&  =\int_{0}^{R}h(r)r^{n+l-1}(u(r;\alpha)^{+})^{p+1}dr.\nonumber
\end{align}
The function $w$ is bounded and eventually monotone, so by Lemma 4.1 of
[K-N-Y] there exists a sequence $\{R_{j}\}$ such that as $j\rightarrow\infty$,%

\begin{equation}
R_{j}\rightarrow\infty,\quad w(R_{j};\alpha)\rightarrow C,\quad R_{j}%
w^{\prime}(R_{j};\alpha)\rightarrow0,\quad R_{j}^{2}w^{\prime\prime}%
(R_{j};\alpha)\rightarrow0. \tag{4.2}%
\end{equation}
Evaluating (4.1) on the sequence $\{R_{j}\}$ we see that the lim sup as
$j\rightarrow\infty$ of the left hand side of (4.1) is $\leq0$ for $p\geq
p_{\ast}$. Evaluation of the right hand side will require several preliminary results.

\begin{lemma}
Assume that $p=p_{\ast}$, and that $(f_{4})$ and $(f_{9})$ are satisfied. Let
$u(r;\alpha)$ be a solution of (1.4). Then for $\gamma>0$ sufficiently small,%
\begin{equation}
lim_{\alpha\rightarrow\infty}\int_{0}^{r_{\alpha}}h(r)r^{n+l-1}(u(r;\alpha
)^{+})^{p+1}dr=\infty. \tag{4.3}%
\end{equation}
\end{lemma}

\begin{proof}
According to Lemma 3.3 $lim_{\alpha\rightarrow\infty}r_{\alpha}=0$, so by
$(f_{4})$ for sufficiently large $\alpha$,%
\begin{align*}
\int_{0}^{r_{\alpha}}h(r)r^{n+l-1}(u(r;\alpha)^{+})^{p+1}dr  &  >\\
\left(  \frac{\alpha}{2}\right)  ^{p+1}\int_{0}^{r_{\alpha}}h(r)r^{n+l-1}dr
&  \geq\\
c\left(  \frac{\alpha}{2}\right)  ^{p+1}\int_{0}^{r_{\alpha}}r^{n+l-1+\gamma
}dr  &  =\\
&  \frac{c}{(n+l+\gamma)2^{p+1}}\alpha^{p+1}r_{\alpha}^{n+l+\gamma}.
\end{align*}
By Lemma 3.3 this is $\geq\alpha^{q}$, for%
\begin{align*}
q  &  =-\{l^{2}+(n+\gamma-\sigma)l+(2\gamma-\sigma n)\}\\
&  =-\{\gamma(2+l)-\sigma(n+l)+l(n+l)\}\\
&  =-\{\gamma(2+l)-(\sigma-l)(n+l)\},
\end{align*}
and $\gamma(2+l)-(\sigma-l)(n+l)<0$ by ($f_{9}$). In fact, it is $<0$ for
$\gamma<\frac{(\sigma-l)(n+l)}{(2+l)}$. Therefore the conclusion of the lemma
follows. Note that a sufficient condition is $\sigma>l$ and $\gamma$
sufficiently small.
\end{proof}

We will also require the following lemma, which gives an apriori bound for
positive solutions on $[r_{0},\infty)$, for some $r_{0}>0$. It is based on
Theorem 3.10 of [N]. Similar results can be found in [B-L], [E-R].

\begin{lemma}
Let $u(r;\alpha)$ be a positive, entire solution of (1.1), with $p\leq
p_{\ast}$, and assume that there exist positive numbers $r_{0},C_{1}$ such
that%
\[
f(r)\geq C_{1}r^{l},\quad-2<l<0,\quad\text{ for }r>r_{0}.
\]
Then there exists a positive constant C, independent of $\alpha$, such that if
$u(r;\alpha)$ is a positive, entire solution then%
\begin{equation}
u(r;\alpha)\leq Cr^{\frac{2-n}{2}},\quad\text{ for }r\geq r_{0},\quad
0<\alpha<\infty. \tag{4.4}%
\end{equation}
\end{lemma}

\begin{proof}
By Green's identity in any ball of radius $R>0$,
\begin{align}
\int_{B_{R}}f\left(  \left|  x\right|  \right)  u\left(  \left|  x\right|
:\alpha\right)  ^{p}dx  &  =-\int_{B_{R}}\Delta u\left(  \left|  x\right|
;\alpha\right)  dx\tag{4.5}\\
&  =\int_{\partial B_{R}}u^{\prime}\left(  r;\alpha\right)  ds\nonumber\\
&  =-\omega_{n}R^{n-1}u^{\prime}\left(  R;\alpha\right)  \text{.}\nonumber
\end{align}
By (4.4) we may assume that%
\[
f(r)\geq C_{1}r^{l},\quad\text{ for }r\geq\frac{r_{0}}{2}.
\]
Then for $R\geq\frac{r_{0}}{2}$, and using the monotonicity of $u(r;\alpha)$,%
\begin{align*}
&  \int_{B_{R}}f\left(  \left|  x\right|  \right)  \left(  u^{+}\left(
\left|  x\right|  ;\alpha\right)  \right)  ^{p}dx\geq\\
&  \int_{B_{R}-B_{\tfrac{r_{0}}{2}}}f\left(  \left|  x\right|  \right)
\left(  u^{+}\left(  \left|  x\right|  ;\alpha\right)  \right)  ^{p}dx\geq\\
&  \int_{B_{R}-B_{\tfrac{r_{0}}{2}}}C_{1}\left|  x\right|  ^{l}\left(
u^{+}\left(  \left|  x\right|  ;\alpha\right)  \right)  ^{p}dx\geq\\
&  C_{1}u\left(  R;\alpha\right)  ^{p}\int_{\tfrac{r_{0}}{2}}^{R}%
r^{n+l-1}dr\geq\\
&  \dfrac{C_{1}}{n+l}u\left(  R;\alpha\right)  ^{p}\left\{  R^{n+l}-\left(
\dfrac{r_{0}}{2}\right)  ^{n+l}\right\}
\end{align*}
Therefore for $R\geq\frac{r_{0}}{2}$,%
\[
\int_{B_{R}}f(|x|)(u^{+}(|x|;\alpha))^{p}dx\geq\frac{C_{1}}{2(n+l)}%
R^{n+l}u(R;\alpha)^{p},
\]
and by (4.5)%
\[
-\omega_{n}R^{n-1}u^{\prime}(R;\alpha)\geq\frac{C_{1}}{2(n+l)}R^{n+l}%
u(R;\alpha)^{p},
\]%
\[
-\frac{u^{\prime}(r;\alpha)}{u(r;\alpha)^{p}}\geq\frac{C_{1}}{\omega_{n}%
}R^{l+1}\quad\text{ for }r\geq\frac{r_{0}}{2}.
\]
Since $l>-2$, an integration from $\frac{r_{0}}{2}$ to $R$ gives%
\begin{align*}
&  -\int_{\tfrac{r_{0}}{2}}^{R}\dfrac{u^{\prime}\left(  r;\alpha\right)
}{u\left(  r;\alpha\right)  ^{p}}dr\\
&  =-\dfrac{u\left(  r;\alpha\right)  ^{1-p}}{1-p}\left.
\raisebox{-0.25ex}[2ex][2ex]{}%
\right|  _{\tfrac{r_{0}}{2}}^{R}\\
&  =\dfrac{1}{p-1}\left\{  u\left(  R;\alpha\right)  ^{1-p}-u\left(
\dfrac{r_{0}}{2};\alpha\right)  ^{1-p}\right\}
\end{align*}
and therefore
\[
\frac{1}{p-1}\left\{  u(r;\alpha)^{1-p}-u\left(  \frac{r_{0}}{2}%
;\alpha\right)  ^{1-p}\right\}  \geq\frac{C_{1}}{\omega_{n}(l+2)}r^{l+2},
\]%
\[
u(r;\alpha)^{1-p}\geq\left(  \frac{C_{1}}{\omega_{n}(l+2)}r^{l+2}+u\left(
\frac{r_{0}}{2};\alpha\right)  ^{1-p}\right)  (p-1)\geq\frac{C_{1}%
(p-1)}{\omega_{n}(l+2)}r^{l+2},
\]
and%
\[
u(r;\alpha)\leq\left(  \frac{C_{1}(p-1)}{\omega_{n}(l+2)}\right)  ^{\frac
{1}{1-p}}r^{\frac{l+2}{1-p}},\text{ for }r\geq r_{0},
\]
and since $p\leq p_{\ast}$, the lemma is proved.
\end{proof}

\begin{proof}
[Proof of Theorem 3]We will evaluate the right hand side of (4.1) using Lemmas
4.1 and 4.2. Note that $r^{0}$ is defined in (3.10).%
\begin{align*}
\left|  \int_{r^{0}}^{R_{j}}h(r)r^{n+l-1}(u(r;\alpha)^{+})^{p+1}dr\right|   &
\leq\\
\int_{r^{0}}^{R_{j}}|h(r)|r^{n+l-1}(u(r;\alpha)^{+})^{p+1}dr  &  \leq\\
\int_{r^{0}}^{\infty}|h(r)|r^{n+l-1}r^{\frac{2-n}{2}(p+1)}dr  &  \leq\\
C\delta\int_{r^{0}}^{\infty}r^{{n+l-1}+\frac{2-n}{2}(p+1)-\beta}  &  \leq\\
C\delta\int_{r^{0}}^{\infty}r^{-(1+\beta)}dr  &  =C_{0},
\end{align*}
by ($f_{6}$) and Lemma 4.2. By the hypotheses $0<\gamma<\gamma_{\ast}$ and
Lemma 4.1, there exists an $\alpha_{1}>0$ such that $\alpha>\alpha
_{1}\Rightarrow\int_{0}^{r^{0}}h(r)r^{n+l-1}(u(r;\alpha)^{+})^{p+1}dr>C_{0}$.
For $\alpha>\alpha_{1}$ by (4.2) we can choose a sequence $\{R_{j}\}$, with
$R_{j}\rightarrow\infty$, such that for $j>j(\alpha)$ the left hand side of
(4.1) is $<C_{0}$, which proves the theorem.
\end{proof}

\begin{proof}
[Proof of Theorem 4]Let
\begin{align*}
h^{+}(r)  &  =max\{0,h(r)\},\\
h^{-}(r)  &  =min\{0,h(r)\}.
\end{align*}
By $(f_{7})$ and $(f_{8})$%
\begin{equation}
0<r_{3}=sup\{r>0:\int_{0}^{r}h(s)s^{n+l-1}ds\leq0\}. \tag{4.6}%
\end{equation}
We choose $k$, $\epsilon$, $\alpha$ and $R_{j}$ as follows:

\begin{enumerate}
\item [(a)]let $0<k<1$, and let $\epsilon>0$ be so small that%
\begin{equation}
\epsilon\int_{r_{0}}^{r_{2}}h^{+}(r)r^{n+l-1}dr<(1-k)\int_{0}^{r_{2}%
}h(r)r^{n+l-1}dr, \tag{4.7i}%
\end{equation}

\item[(b)] let $\alpha$ be so small that%
\begin{equation}
r_{\alpha}>r_{2},\quad\text{ and }(1-\epsilon)^{\frac{1}{p+1}}\alpha\leq
u(r;\alpha)\leq\alpha,\quad\text{ for }0\leq r\leq r_{2}, \tag{4.7ii}%
\end{equation}

\item[(c)] choose $R_{j}>r_{2}$ so large that%
\begin{equation}
L(R_{j},\alpha)<k\alpha^{p+1}\int_{0}^{r_{2}}h(r)r^{n+l-1}dr, \tag{4.7iii}%
\end{equation}
where $L(R_{j},\alpha)$ is the left hand side of (4.1), and assume that
$u(r;\alpha)>0$ for $0\leq r\leq R_{j}$. Then $\int_{0}^{r_{3}}h(r)r^{n+l-1}%
dr=0$, and the right hand side of (4.1) is%
\begin{align}
&  \int_{0}^{R_{j}}h(r)r^{n+l-1}(u^{+}(r;\alpha))^{p+1}dr\tag{4.8}\\
&  =\left\{  \int_{0}^{r_{0}}+\int_{r_{0}}^{r_{3}}+\int_{r_{3}}^{r_{2}}%
+\int_{r_{2}}^{R_{j}}\right\}  h(r)r^{n+l-1}(u^{+}(r;\alpha))^{p+1}%
dr\nonumber\\
&  \geq\alpha^{p+1}\left\{  (1-\epsilon)\int_{0}^{r_{0}}h(r)r^{n+l-1}%
dr+(1-\epsilon)\int_{r_{0}}^{r_{3}}h^{+}(r)r^{n+l-1}dr+\right. \nonumber\\
&  \left.  \int_{r_{0}}^{r_{3}}h^{-}(r)r^{n+l-1}dr+(1-\epsilon)\int_{r_{3}%
}^{r_{2}}h^{+}(r)r^{n+l-1}dr+\int_{r_{3}}^{r_{2}}h^{-}(r)r^{n+l-1}dr\right\}
\nonumber\\
&  \geq\alpha^{p+1}\left\{  (1-\epsilon)\int_{0}^{r_{0}}h(r)r^{n+l-1}%
dr+(1-\epsilon)\int_{r_{0}}^{r_{3}}h^{+}(r)r^{n+l-1}dr\right. \nonumber\\
&  +(1-\epsilon)\int_{r_{0}}^{r_{3}}h^{-}(r)r^{n+l-1}dr+\epsilon\int_{r_{0}%
}^{r_{3}}h^{-}(r)r^{n+l-1}dr+(1-\epsilon)\int_{r_{3}}^{r_{2}}h^{+}%
(r)r^{n+l-1}dr\nonumber\\
&  \left.  +(1-\epsilon)\int_{r_{3}}^{r_{2}}h^{-}(r)r^{n+l-1}dr+\epsilon
\int_{r_{3}}^{r_{2}}h^{-}(r)r^{n+l-1}dr\right\} \nonumber\\
&  =\alpha^{p+1}\left\{  (1-\epsilon)\int_{0}^{r_{2}}h(r)r^{n+l-1}%
dr+\epsilon\int_{r_{0}}^{r_{2}}h^{-}(r)r^{n+l-1}dr\right\} \nonumber\\
&  >k\int_{0}^{r_{2}}h(r)r^{n+l-1}dr,\nonumber
\end{align}
which contradicts (4.7iii).
\end{enumerate}
\end{proof}

\begin{remark}
If the condition $\int_{0}^{r}h(s)s^{n+l-1}ds\geq0$ is satisfied for all
$r>0$, then the conclusion follows from [K-N-Y], Theorem 9.3.
\end{remark}

\section{Proof of Theorem 5}

The proof is based on Theorems 3 and 4. We will construct a specific function
$f(r)$ such that for $p=p_{\ast}$ the corresponding problem (1.4) has
solutions with properties analogous to Example III, Section 1. This result
demonstrates the complexity of the structure of the solution set which can
occur when the hypothesis $r_{H}\leq r_{G}$ of [Y-Y1] is not satisfied.

\begin{lemma}
Under the hypotheses (2.1) the set of slowly decaying solutions of (1.4) is open.
\end{lemma}

\begin{proof}
By Lemma 2.6 of [K-Y-Y] a sufficient condition for the openness of slowly
decaying solutions is that there exist an $r_{2}>0$ such that%
\[
\frac{1}{p+1}r^{l-1}(r^{-l}f(r))_{r}\leq0,\quad r>r_{2}.
\]
This is equivalent to $h(r)\leq0$ for $r>r_{2}$, and is implied by (2.1a).
\end{proof}

Now let $\phi(r;\alpha)$ be the solution of the equation%
\begin{align}
-(r^{n-1}\phi_{r})_{r}  &  =r^{n-1}r^{l}\phi^{p_{\ast}},\tag{5.1}\\
\phi(0)  &  =\alpha,\quad\phi_{r}(0)=0.\nonumber
\end{align}
The solutions $\phi(r;\alpha)$ are all rapidly decaying, and are given
explicitly by%
\begin{equation}
\phi(r;\alpha)=\alpha\left\{  1+\frac{2\alpha^{p_{\ast}-1}}{(p_{\ast
}+1)(n-2)^{2}}r^{(n-2)(p_{\ast}-1)/2}\right\}  ^{-\frac{2}{p_{\ast}-1}}.
\tag{5.2}%
\end{equation}
Consider equation (4.1) with $f\left(  r\right)  $ defined by%
\begin{equation}
f(r)=r^{l}+\epsilon(p_{\ast}+1)r^{l}\int_{0}^{r}s^{-(n+l)}k(s)ds, \tag{5.3}%
\end{equation}
and define $h_{\epsilon}(r)=r(r^{-l}f(r))_{r}=\epsilon(p_{\ast}+1)r^{-(n+l-1)}%
k(r)$. We will compare the solutions of (5.1) with those of equation (1.4)
with $f=f_{\epsilon}$ given by (5.3). \smallskip\noindent By the assumptions
(2.1) the function $f(r)$ satisfies ($f_{4}$), ($f_{6}$), ($f_{7}$) and
($f_{9}$), and clearly $h_{\epsilon}(r)\geq0$ for $r\geq c$. It follows from
Theorems 3 and 4 that for each $\epsilon>0$ there are positive numbers
$\alpha_{0}(\epsilon)$, $\alpha_{1}(\epsilon)$ which satisfy conditions (2.3a)
and (2.3c). We will show that for $\epsilon$ sufficiently small the family of
solutions also satisfies (2.3b), i.e. are slowly decaying. In fact, for
$\alpha_{\ast}>0$ the solution $u(r;\alpha_{\ast})$ converges uniformly on the
interval $[0,c]$ to the solution $\phi(r;\alpha_{\ast})$ as $\epsilon
\rightarrow0$. Therefore there exists an $r_{\ast}>0$ such that for
$r>r_{\ast}$,%
\begin{align*}
P(r;u(r;\alpha_{\ast}))  &  =\frac{1}{p_{\ast}+1}\int_{0}^{r}s^{(n+l-1)}%
h_{\epsilon}(s)(u^{+}(s;\alpha_{\ast}))^{p_{\ast}+1}ds\\
&  \leq\epsilon\left\{  \int_{0}^{a}k(s)(u^{+}(s;\alpha_{\ast}))^{(p_{\ast
}+1)}ds+\int_{a}^{b}k(s)(u^{+}(s;\alpha_{\ast}))^{(p_{\ast}+1)}ds\right. \\
&  \left.  +\int_{b}^{r}k(s)(u^{+}(s;\alpha_{\ast}))^{(p_{\ast}+1)}ds\right\}
\\
&  \leq\epsilon u(0;\alpha_{\ast})^{(p_{\ast}+1)}\int_{0}^{a}k(s)ds+\epsilon
u(b;\alpha_{\ast})^{(p_{\ast}+1)}\int_{a}^{b}k(s)ds\\
&  +\epsilon u(b;\alpha_{\ast})^{(p_{\ast}+1)}\int_{b}^{c}k(s)ds\\
&  =\epsilon\left\{  \phi(0;\alpha_{\ast})^{(p_{\ast}+1)}\int_{0}%
^{a}k(s)ds+\phi(b;\alpha_{\ast})^{(p_{\ast}+1)}\int_{a}^{b}k(s)ds\right. \\
&  \left.  +\phi(b;\alpha_{\ast})^{p_{\ast}+1}\int_{b}^{r}k(s)ds\right\}
+o(\epsilon).
\end{align*}
For $\epsilon$ sufficiently small, by (2.1d) this is $\leq-c^{2}$. The
existence of a slowly decaying solution follows from [Y-Y1, Lemma 2.5]. \smallskip

\noindent From the existence of a slowly decaying solution we can conclude, by
the open property of the slowly decaying solutions under hypotheses 2.1, and
the corresponding property for crossing solutions, that there also exist at
least two rapidly decaying solutions, $u(r;\alpha_{1})$ and $u(r;\alpha_{2})$,
with $\alpha_{1}\in(\alpha_{0},\alpha_{*})$ and $\alpha_{2} \in(\alpha
_{*},\alpha_{3})$.

\begin{remark}
Since $p=p_{\ast}$ the solution $u(r;\alpha_{\ast})$ is slowly decaying by the
Pohozaev identity.
\end{remark}

\end{document}